\documentclass[final]{elsart}
\usepackage{graphicx,amssymb,amsmath,subfig,wrapfig}
\usepackage[round]{natbib}

\makeatletter
\def\enumhook{}

\def\itemhook{}

\def\descripthook{}
\def\enumerate{%
   \ifnum \@enumdepth >\thr@@\@toodeep\else
     \advance\@enumdepth\@ne
     \edef\@enumctr{enum\romannumeral\the\@enumdepth}%
       \expandafter
       \list
         \csname label\@enumctr\endcsname
         {\usecounter\@enumctr\def\makelabel##1{\hss\llap{##1}}%
           \enumhook \csname
enumhook\romannumeral\the\@enumdepth\endcsname}%
   \fi}
\def\itemize{%
   \ifnum \@itemdepth >\thr@@\@toodeep\else
     \advance\@itemdepth\@ne
     \edef\@itemitem{labelitem\romannumeral\the\@itemdepth}%
     \expandafter
     \list
       \csname\@itemitem\endcsname
       {\def\makelabel##1{\hss\llap{##1}}%
         \itemhook \csname
itemhook\romannumeral\the\@itemdepth\endcsname}%
   \fi}

\makeatother

\renewcommand{\itemhook}{\setlength{\topsep}{0pt}
\setlength{\itemsep}{0pt}}
\renewcommand{\enumhook}{\setlength{\topsep}{0pt}
\setlength{\itemsep}{0pt}}

\newtheorem{theorem}{Theorem}
\newtheorem{lemma}[theorem]{Lemma}

\newtheorem{corollary}[theorem]{Corollary}

\newenvironment{proof}{{\bf Proof}}{\hfill\rule{2mm}{2mm}}

\newcounter{algorithmCounter}
\newtheorem{myAlg}[algorithmCounter]{Algorithm}

\def\complaint#1{}
\def\withcomplaints{
\newcounter{mycomplaints}
\def\complaint##1{\refstepcounter{mycomplaints}%
\ifhmode%
\unskip%
{\dimen1=\baselineskip \divide\dimen1 by 2 %
\raise\dimen1\llap{\tiny -\themycomplaints-}}\fi%
\marginpar{\tiny [\themycomplaints]: ##1}}%
}

\newif\ifpdf \ifx\pdfoutput\undefined \pdffalse \else \pdfoutput=1 \pdftrue \fi

\pdfoutput=0

\ifpdf

\else

\fi

\begin{document}
\begin{frontmatter}

\title{Pebble Game Algorithms and Sparse Graphs}

\author[umass]{Audrey Lee\thanksref{nsf1}}
\ead{alee@cs.umass.edu}
\author[smith]{Ileana Streinu\thanksref{nsf2}}
\ead{streinu@cs.smith.edu}

\thanks[nsf1]{Research supported by an
NSF graduate research fellowship and by the NSF grant CCR-0310661 of
the second author.}
\thanks[nsf2]{Research partially funded by NSF grant CCR-0310661.}

\address[umass]{Computer Science Department, University of Massachusetts Amherst, MA, USA}
\address[smith]{Computer Science Department, Smith College, Northampton, MA, USA}

\begin{abstract}
A multi-graph $G$ on $n$ vertices is $(k,\ell)$-sparse if every
subset of $n'\leq n$ vertices spans at most $kn'- \ell$ edges. $G$
is {\em tight} if, in addition, it has exactly $kn - \ell$ edges.
For integer values $k$ and $\ell \in [0, 2k)$, we characterize the
$(k,\ell)$-sparse graphs via a family of simple, elegant and
efficient algorithms called the $(k,\ell)$-pebble games.
%
\end{abstract}


\end{frontmatter}

\section{Introduction}
\label{sec:introduction}
A multi-graph $G=(V,E)$ with $n=|V|$ vertices and $m=|E|$ edges is
{\em $(k,\ell)$-sparse} if every subset of $n'\leq n$ vertices spans
at most $kn'- \ell$ edges.  If, furthermore, $m=kn- \ell$, $G$ is
called {\it tight}. A {\em $(k,\ell)$-spanning} graph is one
containing a tight subgraph that spans the entire vertex set $V$.
For brevity, we will refer to $G$ as a {\em graph} instead of as a
{\em multi-graph} (even though it may have loops and multiple edges)
and will abbreviate {\em $(k,\ell)$-sparse} as {\em sparse}.

\noindent{\bf Historical overview.} Sparse graphs first appeared in
Lor{\'{e}}a \cite{lorea:matroidalFamilies:1979}, as examples of
matroidal families. Classical results of Nash-Williams
\cite{nash-williams:edge-disjoint-spanning-trees:1961} and Tutte
\cite{tutte:decomposing-graph-in-factors:1961} identify the class of
graphs decomposable into $k$ edge-disjoint spanning trees with the
$(k,k)$-tight graphs. Tay
\cite{tay:rigidityMultigraphs-I:1984}
relates them to generic body-and-bar rigidity in arbitrary
dimensions. The $(2,3)$-tight graphs are the {\it generic minimally
rigid} (or Laman) graphs for bar-and-joint frameworks in the plane
\cite{laman:Rigidity:1970}, and the spanning ones correspond to
those which are {\em rigid}.

A $(k,a)$-arborescence is defined as a graph where adding {\em any}
$a$ edges results in $k$ edge-disjoint spanning trees. Results of
Recski \cite{recski:network-I:1984} and Lovasz and Yemini
\cite{lovasz:yemini:genericRigidity:1982} identify Laman graphs with
$(2,1)$-arborescences. For $\ell \in [k, 2k)$, this is extended by
Haas \cite{haas:arboricityGraphs:2002} to an equivalence of
$(k,\ell)$-sparse graphs and $(k,\ell-k)$-arborescences. Whiteley
\cite{whiteley:unionMatroids:1988,whiteley:Matroids:1996} surveys
several rigidity applications where sparse graphs appear, some
having non-integer parameters associated to them. Frank, Szeg{\H{o}}
and Fekete \cite{frank:szego:constructivePackingCoveringTrees:2003,
szego:constructive-sparse:egres-2003,szego:fekete:note-sparse:2005}
study inductive constructions for various subclasses of sparse
graphs, motivated by the so-called Henneberg sequences appearing in
Rigidity Theory \cite{henneberg:graphischeStatik:1911-68}, and Bereg
\cite{bereg:certifying-minimally-rigid-graphs:2005} computes them
with an $O(n^2)$ algorithm for the minimally rigid (Laman) case.

There exist many algorithms for decomposing a graph into
edge-disjoint trees or forests
\cite{edmonds:edge-disjoint-branchings:1973,gabow:westermann:matroidSums:1988,
gabow:matroidApproachPacking-STOC:1991,
tarjan:edgeDisjointSpTrees:1976,tarjan:roskind:edge-disjoint:1985}.
A variation on the $O(n^2)$ time matching-based algorithm of
\cite{hendrickson:uniqueRealizability:1992} for $2$-dimensional
rigidity became the simple and elegant {\it pebble game algorithm}
of Jacobs and Hendrickson \cite{jacobs:hendrickson:PebbleGame:1997},
further analyzed in \cite{berg:jordan:2003}. Practical applications
in studies of protein flexibility led Jacobs et al. \cite{Ja01} to
pebble game heuristics for special cases of three-dimensional
rigidity. However, intriguingly, we have not found anywhere
algorithms applicable to $(k,a)$-arborescences or to the {\em
entire} class of $(k,\ell)$-sparse graphs.

\noindent{\bf Our results.} In this paper, we describe a {\em family
of algorithms}, called the $(k,\ell)$-pebble games, and prove that
they recognize {\em exactly} the $(k,\ell)$-sparse graphs, for the
{\em entire range} $\ell \in [0, 2k)$.


In our terminology, Jacobs and Hendrickson's is a $(2,3)$-pebble
game. We exhibit here the full extent to which their algorithm can
be generalized, and characterize the recognized classes of {\em
graphs}. We study the following fundamental problems.

\begin{enumerate}
\item {\bf Decision:} is $G$ a tight (or just sparse) graph?
\item {\bf Spanning:} does $G$ span a tight subgraph?
\item {\bf Extraction:} extract a maximal sparse subgraph
(ideally, spanning) from a given graph $G$.
\item {\bf Optimization:} from a graph with {\em weighted} edges,
extract a maximum weight sparse subgraph.
\item {\bf Components:} given a non-spanning graph $G$,
find its components (maximal tight induced subgraphs).
%
%
\end{enumerate}

The pebble game algorithms run in time $O(n^2)$ using simple data
structures and induce good algorithmic solutions for all the above
problems. They exhibit the same complexity as Hendrickson's
matching-based algorithm
\cite{hendrickson:thesis:1991,hendrickson:uniqueRealizability:1992}
for $2$-dimensional rigidity. For the {\em special case} of graphs
decomposable into disjoint unions of spanning trees and
pseudo-forests, corresponding to the range $\ell \in [0,k]$ of
$(k,\ell)$-sparse graphs, we remark that there are $O(n^{3/2})$
algorithms due to Gabow and Westerman
\cite{gabow:westermann:matroidSums:1988}. But no better algorithms
than the pebble games are known for the {\it entire} range of
$(k,\ell)$-sparse graphs.

\section{Properties of Sparse Graphs}
\label{sec:sparseProp}
We start by showing why it is {\em natural} to restrict the range of
the integer parameter $\ell$ to $[0,2k)$. Then we identify a dual
property related to a well-known theorem of Nash-Williams and Tutte
\cite{nash-williams:edge-disjoint-spanning-trees:1961,tutte:decomposing-graph-in-factors:1961}
on tree decompositions. Finally, we define {\em components} and give
a detailed characterization of their main structural properties.

All graphs $G=(V,E)$ in this paper have $n=|V|$ vertices and $m=|E|$
edges. For subgraphs $E'\subset E$ induced on subsets $V'\subset V$,
we use $n'=|V'|$ and $m'=|E'|$. The complete multi-graph on $n$
vertices, with multiplicity $a$ on loops and $b$ on edges, is
denoted by $K_n^{a,b}$, and the loopless version by $K_n^b$. The
{\em degree} of a vertex is the number of incident edges, including
loops. The parameters $k$ and $\ell$ are integers.

\paragraph*{Matroidal sparse graphs.}
\label{sec:matroid}
The following Lemma justifies the choice of parameters and points to
a small correction to the informal definition of sparse graphs we
gave in the introduction: because for the range $\ell\in (k, 2k)$
and for $n'=1$, $kn'-\ell$ becomes negative, we should require that
every subset of $n'\leq n$ vertices spans at most $max\{0,kn'-
\ell\}$ edges.

\begin{lemma}\label{thm:sparseProp} {\bf Properties of sparse
graphs.}

\begin{enumerate}
\item \label{lem:upperBoundOnL}
If $\ell \geq 2k$, the class of sparse graphs contains only the
empty graph.
\item \label{lem:lowerBoundOnL}
If $\ell < 0$, the union of two vertex disjoint sparse graphs may
not be sparse.
\item {\bf Loops and parallel edges}
\label{lem:loopsAndParallelEdges}
A sparse graph may contain at most $k-\ell$ loops per vertex. In
particular,  the sparse graphs are loopless when $\ell \geq k$. The
multiplicity of parallel edges is at most $2k- \ell$.
\item {\bf Single vertex graphs}
\label{lem:singleVertex}
In the upper range $\ell\in (k, 2k)$, there are no tight graphs on a
single vertex.
\item {\bf Small tight graphs} {\em (Szeg{\H{o}}
\cite{szego:constructive-sparse:egres-2003})}
\label{lem:smallBlocks}
If $\ell \in [{3\over 2} k, 2k)$ (called the {\em Szeg{\H{o}}
range}), there are no tight graphs on {\em small} sets of $n$
vertices, for $n \in (2, {\ell \over {2k- \ell}})$.
\item {\bf Smallest tight graphs} \label{lem:smallBlocksVal}
When $\ell \in [{3\over 2} k, 2k)$, the smallest non-trivial tight
sparse graphs have $\lceil {\ell \over {2k- \ell}} \rceil$ vertices.
For integer values of ${\ell \over {2k- \ell}}$, there is only one
tight graph on the minimum number of vertices: the complete
multi-graph $K_{\ell \over {2k- \ell}}^{2k- \ell}$; otherwise, there
will be several.
\end{enumerate}
\end{lemma}

\begin{proof}
(1) For $\ell \geq 2k$, any subset of $n'=2$ vertices would span at
most $2k- \ell \leq 0$ edges.
(2) If we take the vertex disjoint union of two tight sparse graphs
on $n_1$, resp. $n_2$ vertices, the union has $n=n_1+n_2$ vertices
and $k(n_1+n_2) - 2l > kn - \ell$ edges, therefore it is not sparse.
(3) Apply the sparsity condition $m' < kn' - \ell$ for $n'=1$ and
$n'=2$.
(4) Indeed, $kn-\ell < 0$ for $n=1$, and the number of edges cannot
be negative.
(5) Assume $\ell\geq k$. A vertex may not span a negative number of
edges, so $n \geq 2$. By part (\ref{lem:loopsAndParallelEdges})
above, a tight graph with $kn- \ell$ edges is a subgraph of the
complete, loopless $(2k- \ell)$-multi-graph $K_n^{2k-\ell}$;
therefore $kn - \ell = m \leq (2k-\ell) {n\choose 2}$. The
inequality between the extremes leads to the condition $f(n) \geq 0$
for the quadratic function $f(n)=a n^2 + b n +c$, with $a= 2k-
\ell$, $b= \ell-4k$ and $c = 2\ell$. The two roots of $f(n)=0$ are
$n_1 = 2 $ and $n_2 = {\ell \over {2k- \ell}}$. The {\em open}
interval between the roots is non-trivial when it contains at least
one integral value, i.e. when $n_2 \geq 3$. This happens exactly
when $\ell > {3\over 2} k$. For values of $n$ within this interval,
all the subgraphs of $K_n^{2k-\ell}$ are $(k,\ell)$-sparse, but none
is tight.
(6) Direct corollary of (5).
\end{proof}

The range of values $\ell\in [0, k)$ is called the {\em lower range}
and $\ell \in [k, 2k)$ is  the {\em upper range}: the threshold case
$\ell=k$ will occasionally be relevant for properties holding in
either range (so we will specify when the lower and upper range
intervals need to be taken as open or closed). The upper range is
further subdivided into two, of which the {\em Szeg{\H{o}} range}
requires special care in applications such as Henneberg sequences.
This phenomenon, of having to deal with special cases depending on
the range of $\ell$, is symptomatic for sparse graphs and impacts
the choice of data structures for our algorithms. At the upper bound
$\ell=2k-1$, the smallest tight graphs are complete graphs. For
example, when $k=3$ and $\ell=5$, the smallest tight graph is $K_5$.
For other values of $k$ and $\ell$, there may be several {\em
smallest} tight graphs. For example, when $k=7$ and $\ell=11$, there
are $6$ smallest tight graphs: all the multi-graphs on $4$ vertices
with a total of $17$ edges and edge-multiplicity at most $3$.

For values of the parameters $k, \ell$ and $n$ in these ranges, we
show now that the tight graphs form the set of {\em bases} of a
matroid. The proof relies on a very simple property of {\em blocks}
given below on page \pageref{sec:blocks}. White and Whiteley, in the
appendix of \cite{whiteley:Matroids:1996}, observed that the matroid
{\em circuit} axioms are satisfied.

\begin{theorem}[{\bf The $(k,\ell)$-sparsity matroid}]
\label{thm:matroid}
Let $n,k$ and $\ell$ satisfy: (1) $\ell \in [0, k]$ and $n \geq 1$;
(2) $\ell \in (k, {3 \over 2} k)$ and $n \geq 2$; (3) $\ell \in [{3
\over 2} k, 2k)$ and $n=2$ or $n \geq {\ell \over {2k- \ell}}$. Then
the collection of all the $(k,\ell)$-tight graphs on $n$ vertices,
is the set of bases of a matroid whose ground set is the set of
edges of the complete multi-graph on $n$ vertices, with loop
multiplicity $k-\ell$ and edge multiplicity $2k-\ell$.
\end{theorem}

\begin{proof}
We verify the three axioms of a basis system. {\bf Equal
cardinality} holds by definition. To prove {\bf Non-emptiness}, we
construct canonical tight graphs as follows. Let $V =
\{1,\ldots,n\}$.  For $\ell \in [0,k)$, $n\geq 1$, place $k- \ell$
loops per vertex; connect the vertices with $\ell$ trees (e.g.
$\ell$ copies of the {\em same} tree). For $\ell \in [k, {3\over 2}
k)$, $n \geq 2$, place $2k- \ell$ parallel edges between vertices
$1$ and $2$. For each vertex $i>2$, place $2k- \ell$ parallel edges
between vertices $i$ and $1$, and $\ell-k < 2k- \ell$ edges between
vertices $i$ and $2$. Finally, consider the case $\ell \in [{3 \over
2} k, 2k)$. For $n=2$, there is only one tight graph, the
$(2k-\ell)$-multi-edge. For $n \geq {\ell \over {2k- \ell}}$, start
with an arbitrary minimum-size tight graph on the set of vertices
indexed from $1$ to $\lceil {\ell \over {2k- \ell}} \rceil$. For all
vertices of larger index $i
> \lceil {\ell
\over {2k- \ell}} \rceil$, place $k$ edges between $i$ and some of
the vertices of index $\leq \lceil {\ell \over {2k- \ell}}
\rceil-1$, saturating the multiplicity $2k- \ell$ of a vertex of
index $i$ before moving on to the next vertex of index $i+1$.

To prove the {\bf Basis exchange axiom}, let $G_j = (V, E_j),
j=1,2$, be two tight graphs and $e_2 \in E_2 \setminus E_1$. We must
show that there exists an edge $e_1 \in E_1 \setminus E_2$ such that
$(V, E_1 \setminus \{e_1\} \cup \{e_2\})$ is tight. Let $e_2 = uv$
(this includes the case $u=v$ when $e_2$ is a loop). Consider all
the tight induced subgraphs (called {\em blocks}) $H_{i} = (V_i,
E_i)$ of $G_1$ containing vertices $u$ and $v$. Let $V' = \bigcap_i
V_i$ and $H' = (V', E')$ be the subgraph of $G_1$ induced on $V'$.
By Theorem \ref{thm:structure}(\ref{lem:blockIntersection}) proved
below in Section \ref{sec:blocks}, $H'$ is a block of $G_1$. Not all
the edges in $H'$ are in $G_2$, i.e. $H'$ cannot be a block of
$G_2$, since $V'$ also spans $e_2$ in $G_2$ and then the subgraph
$E'\cup\{e_2\}\subset E_2$ would violate the sparsity of $G_2$.
Therefore, $H'$ contains at least one edge $e_1\in E_1\setminus
E_2$. We are done if we show that $H_3 = (V', E_1 \setminus \{e_1\}
\cup \{e_2\})$ is sparse. Indeed, $H'$ is the minimal subgraph of
$G_1$ such that the addition of $e_2$ violates sparsity: any other
subset would have been one of the $V_i$, and $V'$ is contained in
it. Since $V'$ is contained in any subset on which sparsity was
violated in $G_1 \cup \{e_2\}$, the removal of $e_1$ restores the
counts.
\end{proof}

In Theorem \ref{thm:matroid}, the ground set $K_n^{k-\ell, 2k-\ell}$
was chosen to produce all the interesting bases. We may enlarge the
ground set, by adding extra loops and parallel edges, or delete
edges from it, by working with a subgraph of $K_n^{k-\ell,
2k-\ell}$, and we still obtain a matroid. In the first case, the
bases will still be restricted to the number of edges required by
the sparsity conditions; in the second case, the bases are maximal
sparse subgraphs of $G$. This allows us later to refer to the {\em
matroidal property of sparse graphs} as reason for the correctness
of the {\em arbitrary order of edge insertion} in the pebble game
algorithms, and of the greedy algorithm for the Optimization Problem
(see \cite{cormen:leiserson:rivest:stein:Algorithms:1995}, p.345 and
\cite{oxley:matroidTheory:1992}).

\paragraph*{Partitioning.}
\label{sec:partition}
Nash-Williams \cite{nash-williams:edge-disjoint-spanning-trees:1961}
and Tutte \cite{tutte:decomposing-graph-in-factors:1961} gave an
alternative definition of $(k,k)$-tight graphs using vertex
partitions and trees: a graph contains $k$ edge-disjoint spanning
trees if and only if every partitioning of the vertex set into $p$
parts has at least $k(p-1)$ edges between them. If, moreover, it has
$kn-k$ edges, it is the edge-disjoint union of $k$ spanning trees
and a $(k,k)$-tight graph. We describe now a slight generalization
of one direction of their criterion, for all $(k,\ell)$-tight
graphs.

\begin{lemma}\label{lem:tightPartition}
Let $G=(V,E)$ be a $(k,\ell)$-tight graph and $P = \{V_1, \ldots,
V_p\}$ a partition of $V$. In the upper range $\ell\in(k, 2k)$,
further assume that each $|V_i|\geq 2$. Then there are at least
$\ell(p-1)$ edges between the partition sets $V_i$.

\end{lemma}

\begin{proof}
Let $E_i$ be the edge set induced by $V_i$ in $G$ and $n_i=|V_i|,
m_i = |E_i|$. By sparsity and the assumption on the size of $V_i$,
$m_i \leq kn_i- \ell, \forall i$ and $\Sigma_i m_i \leq \Sigma_i
(kn_i - \ell) \leq kn-p\ell$. The number of edges between the
partition sets is $m - \Sigma_i m_i \geq kn- \ell-(kn-p\ell) \geq
\ell(p-1)$.
\end{proof}

\begin{lemma}\label{lem:partition2}
Let $G=(V,E)$ be a tight graph. Then every vertex has degree $\geq
k$. Moreover, if $\ell > 0$, then there is at least one edge between
a vertex $v$ and the rest of the vertices $V \setminus \{v\}$.
\end{lemma}

\begin{proof}
If $v \in V$ had degree $d < k$, the induced subgraph on $V
\setminus \{v\}$ would have $kn- \ell-d > kn- \ell - k = k(n-1) -
\ell$ edges, contradicting the sparsity of $G$. This already implies
the second part of the theorem for $\ell\in (k,2k)$, because sparse
graphs in this range have no loops. The other case $ \ell \in (0,
k]$ follows from Lemma \ref{lem:tightPartition}.
\end{proof}

As a simple corollary, when $\ell > 0$, a tight graph is connected.
We will make use of this small observation in Theorem
\ref{thm:structure} (4). Also, as a consequence of the theorem of
Nash-Williams and Tutte we have that, for $\ell \in (0, k]$,
a $(k,\ell)$-tight graph contains $\ell$ edge-disjoint spanning
trees.

\paragraph*{Blocks, Components and Circuits.} \label{sec:blocks}

In a sparse graph, a subset of vertices $V'\subset V$ may span {\it
exactly} $kn'- \ell$ edges, where $n' = |V'|$. In this case, the
induced subgraph is called a {\it block}. A maximal block (with
respect to the set of vertices) is called a {\it component}. We
describe now basic properties of blocks and components.

We start with a decomposition theorem for a sparse graph into {\em
components}, {\em free vertices} (not part of any component) and
{\em free edges} (not spanned by any block, and hence component). In
rigidity applications, the components correspond to rigid clusters.
This decomposition will be used later in speeding up the pebble
game. For this section, denote the range $\ell \in [0,k]$ as the
{\em lower range} and $\ell \in (k,2k)$ as the {\em upper range}.

\begin{figure}[h]
\centering \subfloat[Lower range example ($k=3, \ell=1$): two blocks
overlapping in one vertex.] {\label{fig.blockIntersection31}
\begin{minipage}[b]{0.47\linewidth}
\centering
\includegraphics[scale=.5]{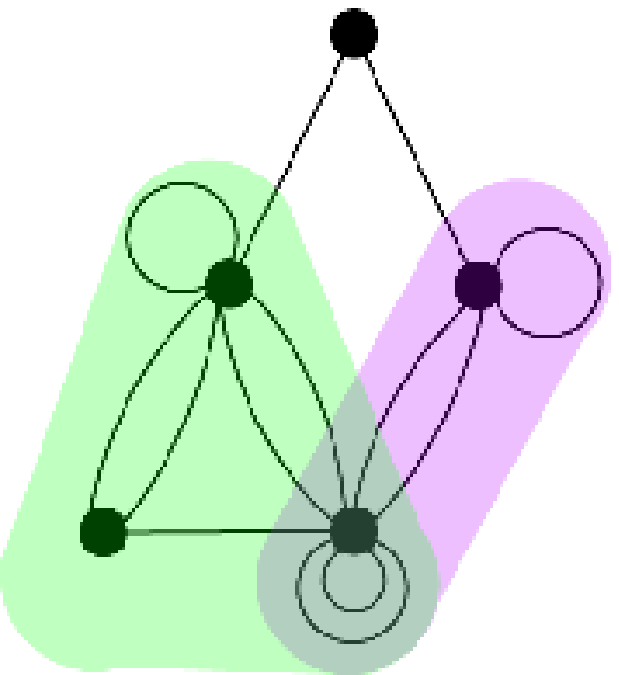}
\end{minipage}}%
\hspace{2mm}
\centering \subfloat[Upper range example ($k=3, \ell=4$): two blocks
overlapping in three vertices.]{\label{fig.blockIntersection34}
\begin{minipage}[b]{0.47\linewidth}
\centering\includegraphics[scale=.5]{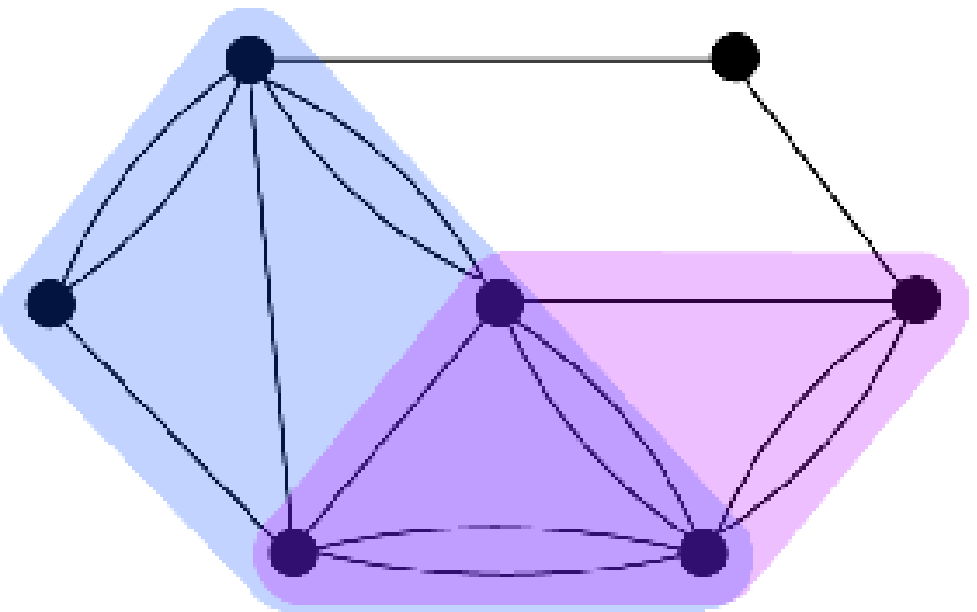}
\end{minipage}}
\caption{Block intersection:  two overlapping blocks whose union and
intersection are also blocks.}
\label{fig.blockIntersection}
\end{figure}

\begin{figure}[t]
\label{fig-decomp}
\centering \subfloat[A $(3,3)$-sparse graph decomposed into
components and free edges.]{\label{fig.ktree3_components_color}
\begin{minipage}[b]{0.22\linewidth}
\centering
\includegraphics[scale=.25]{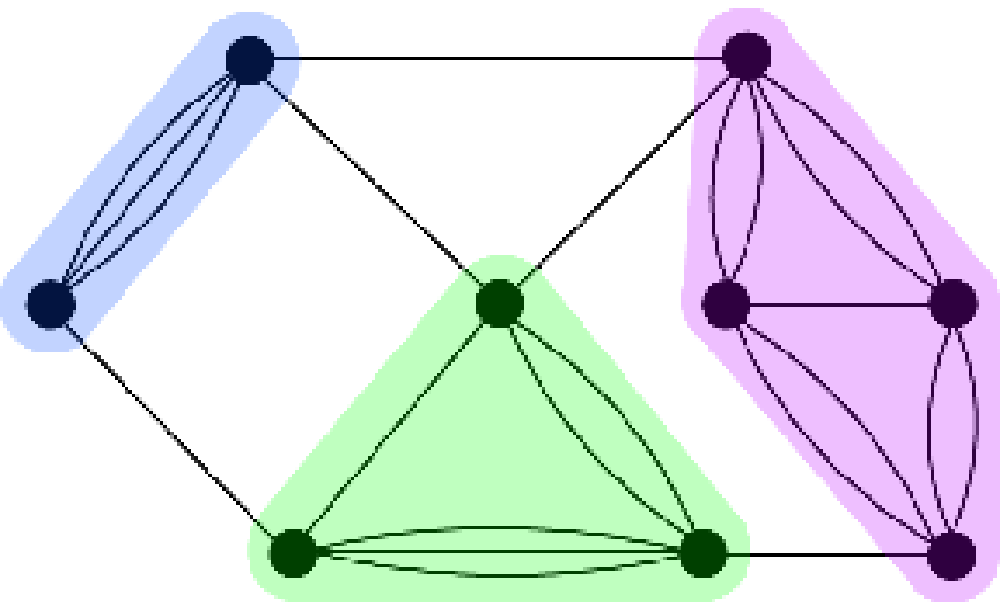}
\end{minipage}}
\hspace{2mm}
\centering \subfloat[A $(2,3)$-sparse graph decomposed into
components.]{\label{fig.laman_components_color}
\begin{minipage}[b]{0.22\linewidth}
\centering
\includegraphics[scale=.25]{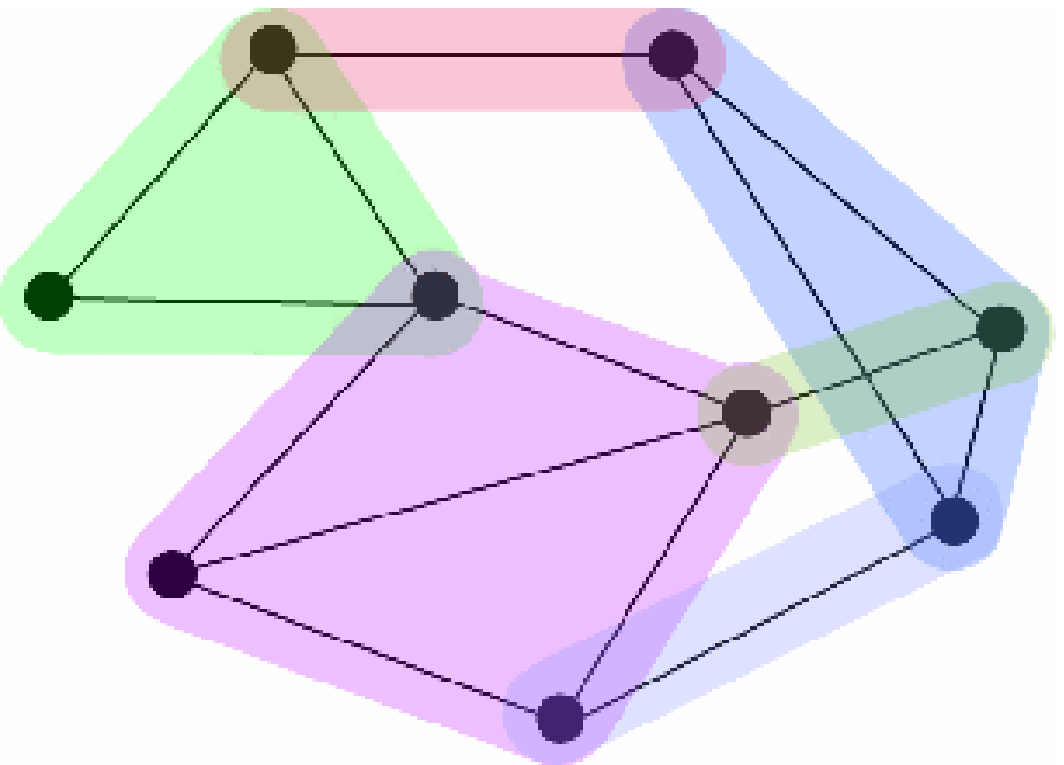}
\end{minipage}}
\hspace{2mm}
\centering \subfloat[A $(2,0)$-sparse graph whose component is not
connected.]{\label{fig.component0}
\begin{minipage}[b]{0.22\linewidth}
\centering
\includegraphics[scale=.25]{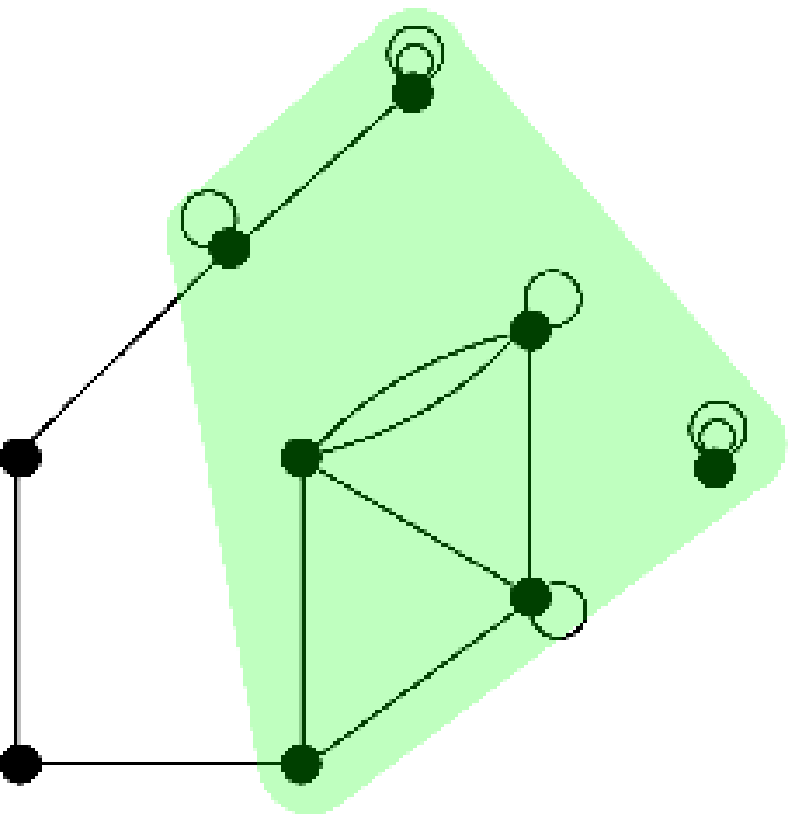}
\end{minipage}}
\hspace{2mm}
\centering \subfloat[A $(3,1)$-sparse graph with $2$ components, a
free vertex and $3$ free edges.]{\label{fig.composition}
\begin{minipage}[b]{0.22\linewidth}
\centering
\includegraphics[scale=.35]{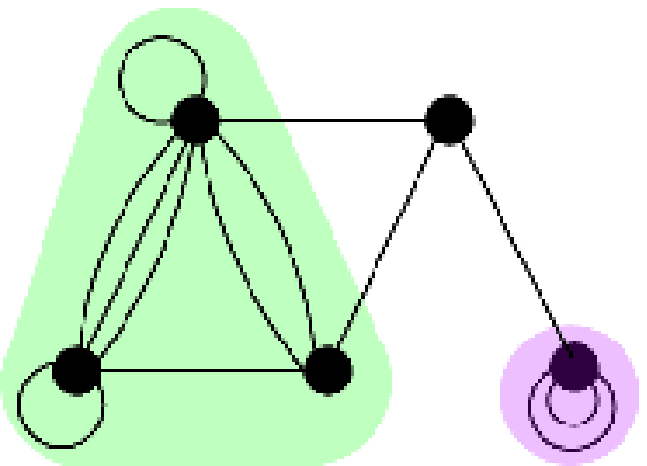}
\end{minipage}}
\caption{Decomposition into components.}
\end{figure}

\begin{theorem}[{\bf Decomposition into Components}]
\label{thm:structure}
Let $G$ be a sparse graph.
\begin{enumerate}
\item \label{lem:blockIntersection} {\bf Block intersection:}
if two blocks intersect in at least: (a) {\em one} vertex, for the
lower range {\em [see Fig. \ref{fig.blockIntersection31}]}; (b) {\em
two} vertices, for the upper range {\em  [see Fig.
\ref{fig.blockIntersection34}]}, then their intersection and union
(with respect to the vertex sets) induce blocks.

\item \label{cor:disjointComp}
{\bf Component interaction:} sparse components are edge-disjoint. In
the lower range, the components are vertex-disjoint {\em [see Fig.
\ref{fig.ktree3_components_color}]}. In the upper range, they
overlap in at most one vertex {\em [see Fig.
\ref{fig.laman_components_color}]}.
\item \label{cor:connected}
{\bf Component connectivity:}
    \begin{enumerate}
    \item \label{lem:comp0}
        When $\ell=0$, there is at most one component, which may not
        be connected
%
        %
 {\em [see Fig.
        \ref{fig.component0}]}.

    \item \label{lem:connectedBlocks}
        When $\ell > 0$, blocks (and therefore components) are connected
%
        %
        {\em [see Fig. \ref{fig.blockIntersection},
        \ref{fig.ktree3_components_color},
        \ref{fig.laman_components_color}, and
        \ref{fig.composition}]}.

    \end{enumerate}

\item {\bf Decomposition:}  $G$ is decomposed
into components, free vertices and free edges.  More specifically:
        \begin{enumerate}
        \item Lower range: a single vertex induces a block
        if and only if it has $k- \ell$ loops. In this case, if
        $\ell = 0$, the block may be a disconnected piece of a larger component,
        otherwise it is a component in itself {\em [see Fig.
        \ref{fig.component0}]}.
        A vertex with fewer than $k- \ell$
        loops in the lower range, or
        a vertex in the upper range is either free or part of a larger block (and hence component).
        When $\ell = 2k-1$, there are no free vertices: each vertex is part of some block (and hence component),
        but it is never a block in itself.
        \item $\ell=k$: a single vertex is loop-free and is always a block. Thus,
        there are no free vertices, and $V$ is partitioned into components
        (possibly connected by free edges){\em [see Fig. \ref{fig.ktree3_components_color}]}.
        \item $\ell=2k-1$: there are no loops or parallel edges. A single vertex is free only when it is
        an isolated vertex of the graph. A single edge is always a block, thus
        there are no free edges, and $E$ is partitioned into
        components {\em [see Fig. \ref{fig.laman_components_color}]}.
        \end{enumerate}
\end{enumerate}

\end{theorem}

\begin{proof}
(1) Let $B_i = (V_i, E_i), i=1,2$, be two blocks of a sparse graph
$G = (V,E)$; they span $m_i = kn_i - \ell$ edges, $i=1,2$. Let
$G_{\cap}$ and $G_{\cup}$ be the subgraphs of $G$ induced on the
intersection $V_1 \cap V_2$ (with $n_{\cap}$ vertices and $m_{\cap}$
edges), resp. union $V_1 \cup V_2$ (with $n_{\cup}$ vertices and
$m_{\cup}$ edges), of their vertex sets. Then $ m_{\cup} = m_1 + m_2
- m_{\cap} = (kn_1 - \ell) + (kn_2 - \ell) - m_{\cap} = k(n_1 + n_2)
- 2\ell - m_{\cap} = k(n_{\cap} + n_{\cup}) - 2\ell - m_{\cap} =
kn_{\cup} - \ell - (m_{\cap} - (kn_{\cap} - \ell))$. Since $G$ is
sparse, $m_{\cup} \leq kn_{\cup} - \ell$; thus, $m_{\cap} -
(kn_{\cap} - \ell) \geq 0$, i.e., $m_{\cap} \geq kn_{\cap} - \ell$.

If $\ell \in [0,k]$, assume $n_{\cap} \geq 1$, and if $\ell \in (k,
2k)$, assume $n_{\cap} \geq 2$. Since $G$ is sparse, $m_{\cap} \leq
kn_{\cap} - \ell$; therefore it follows that $m_{\cap} = kn_{\cap} -
\ell$ and $m_{\cup} = kn_{\cup}- \ell$ and thus, both the induced
intersection and union are blocks.

(2) Follows from the same calculations used in part
(\ref{lem:blockIntersection}).

(3) Lemma \ref{lem:partition2} implies that when  $\ell>0$, tight
graphs are connected. For $(k,0)$-sparse graphs, assume there exist
several vertex-disjoint tight sparse subgraphs (blocks). A simple
application of the sparsity counts shows that the union is also
$(k,0)$-tight.

(4) Take $n=1$ and $n=2$ in the definition of $(k,\ell)$-sparsity,
and analyze each case.
\end{proof}

A reminder that, in matroid theory terminology, a set of elements of
the ground set $E$ of a matroid is {\em independent} if it is a
subset of a basis. An element $e$ of the ground set is {\em
independent with respect to} a given independent set $I\subset E$ if
$I\cup \{e\}$ is an independent set. Thus, sparse graphs are
independent, and independent edges may be added to a sparse graph
until it becomes tight. The obstructions to adding further edges in
a sparse graph are the {\em blocks}, as stated in the following
straightforward corollary to Theorem \ref{thm:matroid}.

\begin{corollary}
\label{cor:indepIFFblock}
An edge is independent with respect to a sparse graph $G$ if and
only if its endpoints do not belong to some block of $G$.
\end{corollary}

A minimal subset (of vertices and edges) violating sparsity is
called a {\em circuit}. An edge which is not independent of a given
sparse graph $G$ violates the sparsity condition on some subset of
vertices and induces a unique circuit, which can be identified using
the criterion below.

\begin{corollary}
\label{cor:minBlockContainingEdge}
Let $G$ be a tight graph and let $e=uv$ be an edge not in $G$. The
intersection of all the blocks containing $u$ and $v$ is a block $H$
of $G$, called the {\em minimal block spanning $e$}. Furthermore,
$H\cup\{e\}$ is a circuit in $G\cup\{e\}$.
\end{corollary}

\section{The basic $(k,\ell)$-Pebble Game Algorithm}
\label{sec:pebbleGames}
\begin{figure}[h]
\centering \subfloat[A {\it Well-constrained} $(3,3)$-pebble game
output, with the final orientation and distribution of the remaining
$3$ pebbles on the input graph.]{\label{fig.3tree}
\begin{minipage}[b]{0.45\linewidth}
\centering
\includegraphics[width=2in]{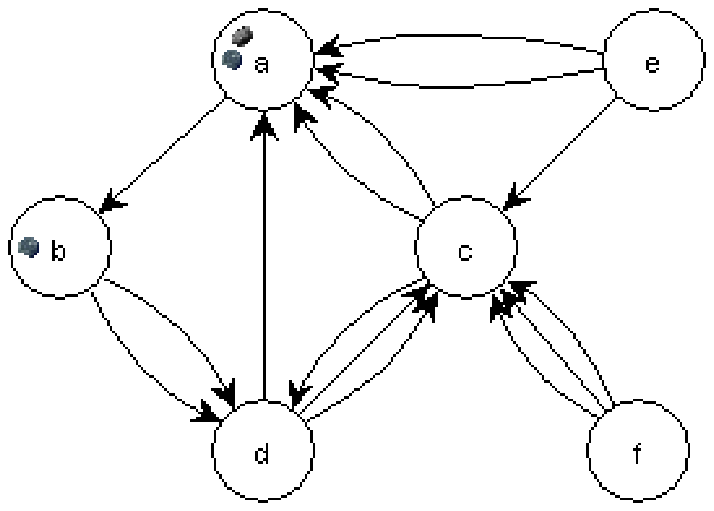}
\end{minipage}}
\hspace{2mm}
\subfloat[An {\it Under-constrained} $(3,3)$-pebble game output:
note the $4$ remaining pebbles. If the dotted edge was part of the
input,  it could not be inserted: the pebble game would {\em
fail}.]{\label{fig.3forest}
\begin{minipage}[b]{0.45\linewidth}
\centering
\includegraphics[width=2in]{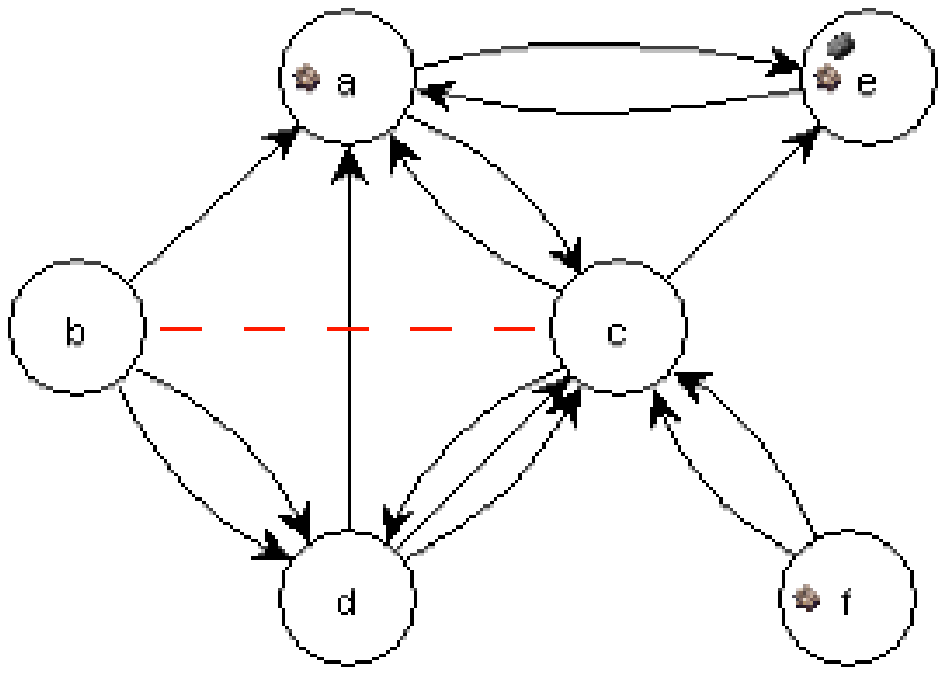}
\end{minipage}}
\caption{Final state of the $(3,3)$-pebble game on two graphs.}
\label{fig.endPebble}
\end{figure}

We turn now to the description of our generalized $(k,\ell)$-{\em
pebble game} for multi-graphs. Fig. \ref{fig.endPebble} illustrates
an example. We start with the simplest version, called the {\em
Basic Pebble Game}. Later, we will extend it to a more efficient
version which takes {\em components} into account. The correctness
of the pebble game as a decision algorithm for sparse graphs is
proven in the next section.

The algorithm depends on two parameters, $k$ and $\ell$: $k$ is the
{\it initial number of pebbles} on each vertex, and $\ell + 1$ is a
lower bound on the total number of pebbles present at the two
endpoints of an edge which is accepted during the execution of the
algorithm.

The algorithm is built on top of a single-person {\em game}, played
on a {\em board} consisting of a set of $n$ nodes, initialized with
$k$ pebbles each. The player inserts edges between the nodes and
orients them. The rules of the game indicate when an edge will be
accepted (and therefore inserted) or rejected, and when the player
can move pebbles and reorient already inserted edges. We give no
rules for when this generic ``game" should be stopped, nor do we
specify what it means to win or to lose it: indeed, we do not
analyze the game per se, but rather the algorithm built on top of
it.

The algorithm takes a given graph as input, and considers its edges
in an arbitrary order. It performs the moves of the game for the
insertion or rejection of each edge. When all the edges have been
considered, the algorithm ends with a classification of the input
graph into one of four categories. The first two, {\it
Well-constrained} and {\it Under-constrained}, correspond to {\it
success} in accepting all the edges of the input graph; the other
two, {\it Over-constrained} and {\em Other}, indicate the failure to
fully accept the input graph. In Section \ref{sec:equivalence} we
prove that these categories correspond exactly to the input graph
being {\em tight}, {\em sparse}, {\em spanning} and {\em neither
sparse nor spanning}. The algorithm is described in
Fig.~\ref{fig:alg:simple}.

\begin{figure}[h]
\noindent \framebox[\textwidth][c]{ 
\parbox{.95\textwidth}
{
\begin{myAlg}{\bf Basic $(k,\ell)$-Pebble Game.}
\label{alg:simple}

\noindent
{\bf Input:} A graph $G=(V,E)$, possibly with loops and multiple edges. \\
{\bf Output:} {\em Well-constrained}, {\em Under-constrained}, {\em
Over-constrained} or {\em Other}.

\noindent{\bf Setup:} Maintain, as an additional data structure, a
directed graph $D$, on which the game is played. Initialize $D$ to
be the empty graph on $V$, and place $k$ pebbles on each vertex. \\

\noindent {\bf Rules:}
\begin{enumerate}
\item {\bf Pebbles.} No more than $k$ pebbles may be present on a vertex at any time.
\item {\bf Edge acceptance.} An edge between two vertices $u$ and $v$
is accepted for insertion in $D$ when a total of at least $\ell + 1$
pebbles are present on the two endpoints $u$ and $v$.
\end{enumerate}

\noindent {\bf Allowable moves:}
\begin{enumerate}
\item {\bf Pebble collection.}
An additional pebble may be collected on a vertex $w$ by searching
the directed graph $D$, e.g., via depth-first search. If a pebble is
found, the edges along the directed path leading to it are reversed
and the pebble is moved along the path until it reaches $w$.
\item {\bf Edge insertion.} If an edge between two vertices $u$ and $v$
is accepted, then at least one of the vertices (say, $u$) contains a
pebble. The edge is inserted in $D$ as a directed edge $u\rightarrow
v$ and a pebble is removed from $u$.
\end{enumerate}
\noindent {\bf Algorithm:} The edges of $G$ are considered in an
arbitrary order and the {\em edge acceptance} condition is checked.
Let $e=uv$ be the current edge. If the acceptance condition is not
met for edge $e$, the algorithm attempts to collect the required
number of pebbles on its two endpoints $u$ and $v$ using the
following strategy: (a) Mark vertices $u$ and $v$ as {\em visited}
for the depth-first-search algorithm (so they will not be searched,
and their pebbles are protected from being moved); (b) Perform
pebble collection using depth-first-search. If one fails to collect
$\ell+1$ pebbles, the edge $e$ is rejected, otherwise it is
accepted.
An accepted edge is immediately inserted into $D$ as specified by
the {\em edge insertion} move.

The algorithm ends when all the edges have been processed. If
exactly $\ell$ pebbles remain in the game at the end, the output is
{\bf Well-constrained} if no edge was rejected, and {\bf
Over-constrained} otherwise. If more than $\ell$ pebbles remain, the
output is {\bf Under-constrained} if there was no edge rejection, or
else {\bf Other}.
\end{myAlg}
}} \caption{Basic $(k,\ell)$-pebble game algorithm.}
\label{fig:alg:simple}
\end{figure}

\noindent{\bf Complexity analysis.} Let $m_a$ be the number of
accepted edges in the final state of the game. Since each accepted
edge requires the removal of one pebble, $m_a = O(kn)$. The only
data structure used by the pebble game is the additional digraph
$D$, whose space complexity is $O(m_a+n) = O(kn)$. Each edge is
considered exactly once and requires at most $\ell+1$ depth-first
searches through $D$, for a total of $O(\ell mn)$ time. For constant
parameters $k$ and $\ell$, and dense input graphs with $O(n^2)$
edges, this algorithm has worst case $O(n^3)$ time and $O(n)$ space
complexity. The time will be improved in Section \ref{sec:pgComp}.

\section{Pebble Game Graphs coincide with
Sparse Graphs}
\label{sec:equivalence}
We are now ready to prove the main theoretical result of the paper,
relating pebble games to sparse graphs.

\begin{theorem}[{\bf Pebble Game Graphs and Sparse Graphs}]
\label{thm:main}
The class of\\
{\em Under-constrained} pebble game graphs coincides with the class
of {\em sparse} graphs, {\em Well-constrained} ones coincide with
{\em tight} graphs, {\em Over-constrained} coincide with {\em
spanning} ones and {\em Other} are neither sparse nor spanning.
\end{theorem}

\begin{corollary}
\label{cor:pbs}
The basic Pebble Game solves the {\em Decision}, {\em Extraction},
{\em Spanning} and, with the slight modification of inserting the
edges in sorted order of their {\em weights}, the {\em Optimization}
problems for sparse graphs.
\end{corollary}

The proof follows from the sequence of lemmas given below. For a
vertex $v$ in the directed graph $D$ at some point in the execution
of the pebble game algorithm, denote by $peb(v)$ the number of free
pebbles on $v$, $span(v)$ the number of loops and by $out(v)$ its
out-degree, i.e. the number of edges starting at $v$ and ending at
{\em a different} vertex (i.e. excluding loops). We extend these
functions to vertex sets in a natural way: for $V'\subset V$,
$peb(V') = \sum_{v\in V'} peb(v)$, $span(V')$ is the number of edges
spanned by $V'$ (including loops) and $out(V')$ is the number of
edges starting at a vertex in $V'$ and ending at a vertex in the
complement $V \setminus V'$.


\begin{lemma} [{\bf Invariants of the Pebble Game}]
\label{lem:invariants}
During the execution of the pebble game algorithm on a graph $G$
with $n$ vertices, for every vertex $v$ and for every subset
$V'\subset V$ on $n'$ vertices, the following invariants are
maintained on $D$. {\em We assume that $n, n'\geq 1$ for $\ell \in
[0, k]$ and $n, n'\geq 2$ for $\ell \in (k, 2k)$.}

\begin{enumerate}
\item \label{kInv} $peb(v) + span(v) + out(v) = k$
\item \label{spanInv}
$peb(V') + span(V') + out(V') = kn'$
\item \label{lInv}
$peb(V') + out(V') \geq \ell$. In particular, there are at least
$\ell$ free pebbles in the digraph $D$.
\item \label{sparseInv} $span(V') \leq kn' - \ell$
\end{enumerate}
\end{lemma}

\begin{proof}
(1) The invariant obviously holds when the game starts. When an edge
is inserted into $D$ and is oriented away from $v$, a pebble is
removed from $v$; this is true for loops as well, so the total sum
is maintained. During a pebble search, if $v$ lies along a path that
is reversed to bring a pebble to the path's source, $out(v)$ remains
unchanged. If $v$ is the source of a path reversal, $out(v)$ is
decreased by $1$ and $peb(v)$ is increased by $1$; if $v$ is the
target of a path reversal, $out(v)$ is increased by $1$ and $peb(v)$
is decreased by $1$. Hence the sum $peb(v) + span(v) + out(v)$
remains constant throughout the game.

(2) If $m_1$ is the number of non-loop edges spanned by $V'$, then
$out(V') = \Sigma_{v \in V'} out(v) - m_1$ and $span(V') = m_1 +
\Sigma_{v \in V'}span(v)$. Therefore $peb(V') + span(V') + out(V') =
peb(V') + (m_1 + \Sigma_{v \in V'}span(v)) + ((\Sigma_{v \in V'}
out(v)) - m_1) = \Sigma_{v \in V'} peb(v) + \Sigma_{v \in V'}span(v)
+ \Sigma_{v \in V'} out(v)  = \Sigma_{v \in V'} (peb(v) + span(v) +
out(v)) = kn'$ (by Invariant (\ref{kInv}).

(3) When the game starts, there are no outgoing edges from $V'$ and
$peb(V')=kn'\geq \ell$. Consider now the last time an edge incident
to $V'$ was inserted or reoriented by the pebble game algorithm.
Four cases have to be analyzed: if the edge had both endpoints in
$V'$, if it went between $V'$ and $V\setminus V'$ and oriented away
from or towards $V'$, or if it was an edge reorientation. In the
first case, at least $\ell$ pebbles must be present on the endpoints
of the edge after the insertion, so $peb(V') \geq \ell$. In the
second case, the invariant was true before the insertion, and it is
true after the insertion: if the edge was inserted away from $V'$, a
pebble was consumed from $V'$ but an outgoing edge was inserted; in
the other case, the number of pebbles and outgoing edges was not
modified. Finally, if the last move was an edge reorientation, it
was either bringing in a pebble from the outside to the inside of
$V'$, and decreasing by $1$ the number of outgoing edges, or
vice-versa, when it was decreasing the number of inside pebbles by
one and increasing the number of outgoing edges.

(4) Straightforward, since $span(V') = kn'- (peb(V') + out(V'))$,
and $peb(V') + out(V') \geq \ell$.
\end{proof}

The following corollaries follow directly from Invariant
\ref{sparseInv}.
\begin{corollary}
For any subset $V' \subseteq V$, $V'$ spans a block if and only if
$peb(V') + out(V') = \ell$.
\end{corollary}

\begin{corollary}
\label{cor:pgGraphsAreSparse}
{\em Under-constrained} pebble game graphs are {\em sparse}, {\em
Well-constrained} ones are {\em tight}, {\em Over-constrained} ones
are {\em spanning}.
\end{corollary}

This completes the proof of one direction, characterizing the
sparsity of the graphs classified by the algorithm. We move now to
prove the other direction, that the algorithm classifies correctly
sparse, tight and spanning graphs. Denote by $Reach(v)$ the {\it
reachability region} of a vertex $v$ (at some point during the
execution of the algorithm): the set of vertices that can be reached
via directed paths from $v$ in $D$. For example, in Figure
\ref{fig.3forest}, $Reach( d ) = \{a, c, d, e\}$.


\begin{lemma}
\label{lem:indepIncreasePeb}
If $e = uv$ is independent (but not yet inserted) in $D$, and
strictly fewer than $\ell+1$ pebbles are present on $u$ and $v$, a
pebble can be brought to one of $u$ or $v$ without changing the
pebble count of the other vertex.
\end{lemma}

\begin{proof}
Let $V' = Reach(u) \cup Reach(v)$; $e$ is independent, so $span(V')
< k|V'|- \ell$. Since $V'$ is a union of reachability regions,
$out(V') = 0$. By Lemma \ref{lem:invariants}, Invariant \ref{lInv},
$peb(V') > \ell$. By assumption, $peb(u) + peb(v) < \ell+1$. Then
there exists $w \in V'$ such that $w \not= u$ and $w \not=v$ with at
least one free pebble. If $w \in Reach(u)$, bring the pebble from
$w$ to $u$. Otherwise, $w \in Reach(v)$; bring the pebble from $w$
to $v$.
\end{proof}

\begin{lemma}
\label{lem:indepInsert}
An edge is inserted by the pebble game if and only if it is
independent in $D$.
\end{lemma}

\begin{proof}
Let $e = uv$ be an edge of $G$, not yet inserted into $D$, the
current state of the directed graph pebble game data structure. By
applying Lemma \ref{lem:indepIncreasePeb} repeatedly, it follows
that $\ell+1$ pebbles can be gathered on the endpoints $u$ and $v$;
thus, $e$ will be inserted into $D$ by the pebble game.
\end{proof}

%
\begin{figure}[h]
\centering 
\includegraphics[width=1.in]{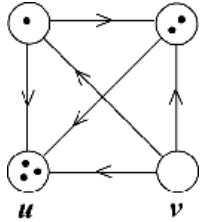}
\vspace{-0.1in}
\caption{A $(3,5)$-pebble game where no edge parallel to $uv$ can be
inserted.}
\label{fig:counterEx} \end{figure}
%
%
It is instructive to notice that it does not suffice to require that
$\ell+1$ pebbles be present in the reachability regions of $u$ and
$v$. In Fig. \ref{fig:counterEx}, an example of a $(3,5)$-pebble
game is shown. The reachability region for the pair $u$ and $v$
contains $6$ pebbles, but not all can be collected on the two
vertices $u$ and $v$. No edge parallel to $uv$ can be inserted. Note
also that the reachability region of a vertex may change after a
pebble move; the previous proof requires the independence of $e$ at
each application of Lemma \ref{lem:indepIncreasePeb}.

\begin{lemma}
\label{lem:sparseArePG}
The pebble game returns {\em Under-constrained} for sparse, but not
tight graphs, {\em Well-constrained} for tight ones, {\em
Over-constrained} for spanning graphs and {\em Other} for graphs
that are neither spanning nor sparse.
\end{lemma}

\begin{proof}
Let $G$ be a sparse graph with $n$ vertices and $m$ edges. Because
sparse graphs form a matroid (Theorem \ref{thm:matroid}), the order
in which the edges are considered can be arbitrary. By Lemma
\ref{lem:indepInsert}, every independent edge is inserted by the
pebble game. Thus, the pebble game is successful on sparse graphs.

If $G$ is sparse, but not tight, $m < kn- \ell$. By Lemma
\ref{lem:invariants}, Invariant \ref{lInv}, the number of free
pebbles in the final game graph must be $> \ell$ and the result is
{\it Under-constrained}. If $G$ is tight, $m = kn- \ell$ and the
number of free pebbles is exactly $\ell$; the game returns {\it
Well-constrained}. If $G$ is spanning, it contains a tight subgraph
which will be accepted, after which there won't be enough pebbles
and the remaining edges will be rejected; the result is {\em
Over-constrained} in this case. If $G$ is neither spanning nor
sparse, there must be $> \ell$ pebbles in the final game as well as
at least one dependent (and thus rejected) edge.
\end{proof}

Corollary \ref{cor:pgGraphsAreSparse} and Lemma
\ref{lem:sparseArePG} prove Theorem \ref{thm:main}.

\section{Component Pebble Games}
\label{sec:pgComp}
The graph $D$ maintained by the basic pebble game algorithm is
sparse. It can therefore be decomposed into components. We now
present a modification of the basic pebble game that maintains and
uses these components to obtain an algorithm one order of magnitude
faster.

\noindent{\bf Pebble Game with Components.} Its input, output and
additional directed graph $D$ are the same as for the basic Pebble
Game. We give first the overall structure of the algorithm in Figure
\ref{fig:alg:efficient}; additional subroutines and some
implementation details will be described next.

\begin{figure}
[h]\noindent \framebox[\textwidth][c] {
\parbox{.95\textwidth}
{
\begin{myAlg}{\bf Component Pebble Game.}
\label{alg:efficient}

\noindent
{\bf Input:} A graph $G=(V,E)$, possibly with loops and multiple edges. \\
{\bf Output:} {\em Well-constrained}, {\em Under-constrained},
{\em Over-constrained} or {\em Other}. \\
{\bf Method:} Play the basic pebble game with the following
modifications. Maintain components throughout the game. When
considering edge $e=uv$, first check if $u$ and $v$ are in some
common component or if $u=v$ (i.e., $e$ is a loop) and $\ell \in [k,
2k)$. If so, {\bf reject} and discard $e$. Otherwise, perform pebble
searches to gather $\ell+1$ pebbles on $u$ and $v$, and insert edge
$e$. Detect new component, if one is formed, and perform necessary
component maintenance.
\end{myAlg}
}} \caption{Component pebble game algorithm.}
\label{fig:alg:efficient}
\end{figure}

If the endpoints of an edge do not belong to a component, the pebble
searches are guaranteed to succeed, so the edge will be accepted. A
newly inserted edge may be a free edge or lead to the creation of a
new component (possibly by merging some already existing ones). We
present next two different algorithms for computing the vertex set
of the new component. Algorithm \ref{alg:detectComp1}, shown in
Figure \ref{fig:alg:detectComp1}, generalizes the approach of
\cite{jacobs:hendrickson:PebbleGame:1997} and works similarly to
breadth-first search on {\em in-coming edges} of an already detected
block. An example is shown in Figure \ref{fig.componentDetection33};
notice that vertex $f$, although it has an edge directed towards
$Reach(a,c)$, contains a pebble and is not added to the component.

\begin{figure}[h]
\noindent \framebox[\textwidth][c] {
\parbox{.95\textwidth}
{
\begin{myAlg}
\label{alg:detectComp1}
{\bf Component detection I}

\noindent {\bf Input:} Directed pebble game graph $D = (V,E')$, into
which edge $e = uv$ has  just been inserted. At least $\ell$ pebbles
are present on $u$ and $v$. If $\ell = 0$, the vertex set $V_0$
(which may be empty) of the single component of $D\setminus \{e\}$
is also given.\\
{\bf Output:} The vertex set $V'$ of the new component induced by
$e$, if
one was formed; $\emptyset$, otherwise.\\
{\bf Method:}
\begin{enumerate}
\item If more than $\ell$ pebbles are present on $u$ and $v$,
{\bf return $\emptyset$}: the new edge is {\em free}.
\item Otherwise, compute $Reach(u,v) = Reach(u) \cup Reach(v)$.
    \begin{enumerate}
    \item If any $w \in Reach(u,v)$ has at least one free pebble,
    {\bf return $\emptyset$}.
    \item Otherwise, initialize $V' = Reach(u,v)$.
        Initialize queue $Q$, and enqueue all vertices in
        $V\setminus V'$ with an edge {\em into} $V'$. \\
        While $Q$ has elements
        \begin{enumerate}
        \item Dequeue vertex $w$ from $Q$.
        \item Compute $Reach(w)$.
        \item If all vertices in $Reach(w)$ (other than $u$ and $v$) have
        no free pebbles
            \begin{enumerate}
            \item Set $V' = V' \cup Reach(w)$.
            \item Enqueue all vertices (that have not been previously
            enqueued) with an edge into $Reach(w)$.
            \end{enumerate}
        \end{enumerate}
    \end{enumerate}
\item If $\ell = 0$, merge into $V'$ the vertices of the existing component of $G$
(if it exists).
\item {\bf Return $V'$.}
\end{enumerate}
\end{myAlg}
}} \caption{Component detection algorithm I.}
\label{fig:alg:detectComp1}
\end{figure}

%
\begin{figure}[h]
\centering \subfloat[Edge $ac$ is successfully
inserted.]{\label{fig.componentDetection33_0}
\begin{minipage}[b]{0.3\linewidth}
\centering
\includegraphics[width=1.6in]{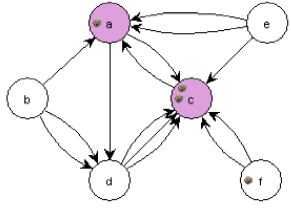}
\end{minipage}}
\hspace{2mm}
\subfloat[$Reach(a,c)$ is detected to have exactly 3
pebbles.]{\label{fig.componentDetection33_1}
\begin{minipage}[b]{0.3\linewidth}
\centering
\includegraphics[width=1.6in]{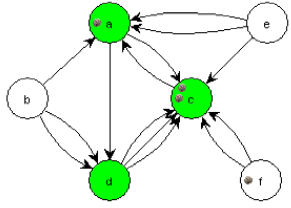}
\end{minipage}}
\hspace{2mm}
\subfloat[Component $\{a,b,c,d,e\}$ is
detected.]{\label{fig.componentDetection33_2}
\begin{minipage}[b]{0.3\linewidth}
\centering
\includegraphics[width=1.6in]{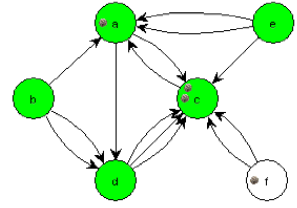}
\end{minipage}}
\caption{$(3,3)$ component detection after edge $ac$ has just been
inserted. First, $Reach(a,c) = \{a,c,d\}$ is detected as a block as
it contains exactly 3 pebbles; then, component $\{a, b, c, d, e\}$
is detected.}
\label{fig.componentDetection33}
\end{figure}

Figure \ref{fig:alg:detectComp2} describes the second component
detection algorithm, which generalizes \cite{berg:jordan:2003}. It
works by finding the {\em complement} of the vertex set of the newly
formed component, and does not require special treatment for
$\ell=0$.

\begin{figure}[h]
\noindent \framebox[\textwidth][c]{
\parbox{.95\textwidth}
{
\begin{myAlg}
\label{alg:detectComp2}
{\bf Component detection II}

\noindent {\bf Input:} Directed pebble game graph $D = (V,E')$, into
which edge $e = uv$ has  just been inserted. At least $\ell$ pebbles
are present on $u$ and $v$.\\
{\bf Output:} The vertex set $V'$ of the new component induced by
$e$, if one was formed; $\emptyset$, otherwise.\\
{\bf Method:}
\begin{enumerate}
\item If more than $\ell$ pebbles are present on $u$ and $v$,
{\bf return $\emptyset$}: the new edge is {\em free}.
\item Otherwise, compute $Reach(u,v) = Reach(u) \cup Reach(v)$.
    \begin{enumerate}
    \item If any $w \in Reach(u,v)$ has at least one free pebble,
    {\bf return $\emptyset$}.
    \item Otherwise, let $D'$ be the directed graph obtained from $D$ by
    reversing the direction of every edge. For all vertices
    $w \in V \setminus Reach(u,v)$ with at least one free pebble,
    perform a depth-first search in $D'$ from $w$. {\bf Return $V'$},
    the set of non-visited vertices from all these searches.
    \end{enumerate}
\end{enumerate}
\end{myAlg}
}} \caption{Component detection algorithm II.}
\label{fig:alg:detectComp2}
\end{figure}

\noindent{\bf Component maintenance.} Maintaining components
requires additional bookkeeping. By Theorem \ref{thm:structure}, we
must take the range of $\ell$ into account. When $\ell = 0$, there
is at most one component, which is maintained by a simple {\em
marking scheme:} a vertex is marked if and only if it lies in the
component. When $\ell \in (0, k]$, the components are vertex
disjoint. Their maintenance is accomplished with a simple {\em
labeling scheme:} each vertex is labeled with an id of the component
to which it belongs. In the upper range, when $\ell\in (k,2k)$,
components may overlap in a single vertex. We maintain a list of the
components, represented by their vertex sets, as well as an $n
\times n$ matrix. The matrix is used to provide constant time
queries for whether two vertices belong to some common component;
there is a $1$ in entry $[i,j]$ if such a component exists and a $0$
otherwise.

When a new component on $V'$ has been detected, we must perform the
necessary bookkeeping to update the data structures. When $\ell \in
[0,k]$ (the {\em lower range}), we simply update the marks or labels
of vertices in $V'$ to record the newly detected component. For the
{\em upper range}, we first mark all vertices in $V'$. Then, for
each previous component $V_i$, all of whose vertices have been
marked, delete $V_i$ from the list of components and update the
matrix.


We are now ready to state and prove that component pebble games
correctly solve most of the fundamental problems presented in the
Introduction: {\em Decision}, {\em Spanning}, {\em Extraction}, {\em
Optimization} and {\em Components}.

\begin{theorem}
\label{thm:efficient}
The graphs recognized by component pebble games are the same as
graphs recognized by basic pebble games, and components are
correctly computed.
\end{theorem}

\begin{proof}
The component pebble game differs from the basic pebble game by
maintaining components and rejecting edges precisely when both
endpoints lie in a component. Thus, by Corollary
\ref{cor:indepIFFblock}, the component pebble game accepts an edge
if and only if it is independent.

To show that the component pebble games correctly maintain
components, observe that Algorithm \ref{alg:detectComp1} detects a
maximal connected subgraph (with respect to the vertices) with no
outgoing edges in which exactly $\ell$ pebbles are present. By Lemma
\ref{lem:invariants}, Invariant \ref{spanInv}, this subgraph must be
a block. When $\ell > 0$, by Theorem
\ref{thm:structure}(\ref{lem:connectedBlocks}), components are
connected; thus, Algorithm \ref{alg:detectComp1} detects a
component. When $\ell=0$, there may be at most one component by
Theorem \ref{thm:structure}(\ref{lem:comp0}) since the union of two
blocks is a block, and the algorithm computes it.

The visited vertices in Algorithm \ref{alg:detectComp2} are those
that can reach a pebble (in the original orientation) on a vertex
other than $u$ and $v$, and thus form the complement of the unique
component containing $e$.
\end{proof}

The time complexity on the algorithm is now $O(n^2)$ as dependent
edges are rejected in constant time. Component detection and
resulting updating of the data structures can be accomplished in
linear time. More specific, implementation-related details on how to
actually achieve this for the {\em upper range} are given in
\cite{streinu:lee:theran:findingRigidComponents:2005} using a
similar data structure called {\em union pair-find}; while {\em
union pair-find} maintains edge sets, the pebble game algorithm does
not need to do it explicitly, and thus the associated implementation
details for edge sets can be ignored.

Space complexity is linear for the {\em lower range} and $O(n^2)$ in
the {\em upper range}, due to its additional matrix. An alternative
solution to {\em union pair-find} presented in
\cite{streinu:lee:theran:findingRigidComponents:2005} uses only
$O(n)$ space, though it requires the edges to be considered in a
specific order and does not solve the Optimization problem.

\section{Applications}
\label{sec:henn}

\noindent{\bf Henneberg Sequences.} Originating in
\cite{henneberg:graphischeStatik:1911-68} (see also
\cite{tay:whiteley:generatingIsostatic:1985}), these are inductive
constructions for Laman graphs and other classes of rigid
structures. We extend the concept to tight graphs: at the base case,
start with a small tight graph; each inductive step would create a
tight graph with an additional vertex by specifying $b$ edges for
removal before adding the new vertex of degree $k+b$. In addition,
$b$ can be chosen to be small $b\in[0,k]$.

We remind the reader of the {\bf matroidal conditions} on tight
graphs: (1) $\ell \in [0, k]$ and $n \geq 1$, or (2) $\ell \in (k,
{3 \over 2} k)$ and $n \geq 2$, or (3) $\ell \in [{3 \over 2} k,
2k)$ and $n \geq \lceil {\ell \over {2k- \ell}} \rceil.$ We refer to
the smallest values of $n$ as the {\bf base-case conditions}; when
$n$ is strictly larger, we call them the {\bf non-triviality
conditions}. The following lemma is the key to proving the existence
of a {\em Henneberg reduction}: given a tight graph, remove vertices
one at a time until a base case is reached. This leads to an
quadratic algorithm for computing the entire sequence.

\begin{lemma} \label{lem:flipNeighbor}
Let $v$ be a vertex of degree $k+b > k$ in a tight graph $G$. Then,
after the removal of any edge $e = uv$, there exists a new edge
whose insertion results in a tight graph. If $\ell \in [0, {3 \over
2} k)$, this edge can be found among the neighbors of $v$;
otherwise, it is found in a larger set containing the neighbors of
$v$ whose size satisfies the {\bf base-case conditions}.
\end{lemma}

\begin{proof}
Consider the sparse graph after the removal of $e$; it is broken
into components, free edges and free vertices. Let $V'$ be the
neighbors of $v$ (but not $v$ itself). If $\ell \in [{3 \over 2},
2k)$, add enough vertices to $V'$ to satisfy the {\bf base-case
conditions}.

We claim that the vertex set $V'$ cannot form a block; in fact, it
cannot span more than $k|V'|-\ell - b$ edges. Indeed, suppose, for a
contradiction, that $V'$ spanned more than $k|V'| - \ell -b$ edges.
Since the degree of $v$ is $k+b$, the size of the induced set of
edges in $G$ on $V' \cup\{ v\}$ is more than $k|V'|-\ell-b + k + b =
(k|V'| - \ell) + k = k(|V'| + 1) - \ell = k|V'\cup\{v\}| - \ell$.
This contradicts the sparsity of $G$.

Since $V'$ does not form a block and its number of vertices
satisfies the {\bf base-case conditions}, it is not saturated with
edges. Therefore, because of the matroidal property of base
extension, there exists an edge not already spanned by $V'$, which
can be added to restore tightness.
\end{proof}

It is a simple exercise to show the existence of a vertex with
bounded degree in $[k,2k]$: indeed, the average degree in a sparse
graph is at most $2k$, and each vertex $v$ has degree at least $k$
(or else sparsity would be violated on $V\setminus\{v\}$). We can
then apply Lemma \ref{lem:flipNeighbor} $O(k)$ times repeatedly to
compute a single Henneberg reduction step.
The Henneberg sequence is obtaining by iterating the Henneberg
reduction step until we reach a base case.

This leads directly to an $O(n^2)$ algorithm for solving the
Henneberg reduction problem by using the pebble game. Figure
\ref{fig:alg:henneberg} describes one step of the algorithm. Each
edge removal puts back one pebble and searches in a constant-size
vertex subset for at most $O(k^2)$ possibilities of edge-insertion,
taking a total of $O(n)$ time in the necessary pebble searches.

\begin{figure}[h]
\noindent \framebox[\textwidth][c]{
\parbox{.97\textwidth}
{
\begin{myAlg}{\bf Henneberg reduction step.}
\label{alg:henneberg}\\
{\bf Input:} The directed graph produced by the pebble game, played
on a tight graph $G$ satisfying the {\bf non-triviality
conditions}.\\
{\bf Output:} A Henneberg reduction step for $G$.\\
{\bf Method:} Find a vertex of degree $k+b$, with $ b\in [k, 2k)$.
If $\ell = 0$, $b$ may also be $k$. Compute the neighbor set $V_v$
of $v$. If $\ell \in [{3 \over 2}k, 2k)$, let $V'$ be {\em any} set
of size $\lceil {\ell \over {2k- \ell}} \rceil$ that includes $V_v$,
else $V'=V_v$.

Repeat $b$ times:\\
Use the pebble game to find an edge with endpoints in $V'$ which is
not already spanned by $V'$ and not in a component. Insert it.
\end{myAlg}}}
\caption{Henneberg reduction step algorithm.}
\label{fig:alg:henneberg}
\end{figure}

\noindent{\bf Circuits and redundancy.} A graph $G$ is said to be
{\em $(k,\ell)$-redundant} if it is spanning and the removal of {\em
any} edge produces a graph which is still spanning. A {\em circuit}
is a special type of redundant graph, where the removal of any edge
produces a tight graph.

We can detect a circuit associated with a dependent edge $e=uv$ with
respect to $D$ during the pebble game by collecting $\ell$ pebbles
on $u$ and $v$ and computing $Reach(u,v)$, which is done in linear
time; the edges in $D$ spanned by $Reach(u,v)$ along with $e$
comprise the circuit.

To decide redundancy of the input graph $G$, simply detect circuits
during the game and mark all the edges in circuits, as they are
computed; if all edges are marked at the end of the game, the graph
is redundant. If the graph is not redundant, unmarked edges are {\em
bridges}; after their removal, the vertex sets of the sparsity
components in the resulting graph correspond to the vertex sets of
{\em redundant components}: induced subgraphs that are redundant.
These algorithms run in $O(mn)$ time.


\setlength{\bibsep}{.04em}
\renewcommand{\baselinestretch}{0}

\label{sec:biblio}

\end{document}